\documentclass[11pt,draft]{article}
\usepackage[cp1250]{inputenc}
\usepackage{latexsym,amssymb,amsmath}

\title{On some criterions of finiteness conditions in commutative
Moufang loops}
\date{23.04.2008}
\author{A. D. Babiy, N. I. Sandu}
\begin{document}
\maketitle

\begin{abstract}

The various finiteness conditions in  commutative Moufang loops
are characterized using the notions of centra\-lizer of subloops
and cen\-tra\-li\-zer of subgroups of its multiplication group.
\smallskip\\
 \end{abstract}

The (\textit{special}) \textit{rank} of loop $L$
 is called the least positive number $rL$ with the following
 feature: any finitely generated subloop of loop $L$ can be
 generated by $rL$ elements; if there are not such numbers, then
 we suppose that $rL = \infty$. Let $\Omega$ denote one of following
classes of loops: the class of finite loops; the class of finitely
generated loops; the class of loops of finite rank; the class of
loops with maximum condition for its subloops;  the class of loops
with minimum condition for its subloops.

It is known that in some classes of groups, loops   various
finiteness condi\-tions of its centre are transferred on these
groups, loops. For example, in [1] it is proved that if the centre
of finitely generated nilpotent group is finite, then the group
itself is finite. Further, for a commutative Moufang $ZA$- loop
$L$ with multiplication group $\frak M$ from paper [2 - 4] it is
follows the equity of statements: 1)  $L$ belongs to class
$\Omega$; 2) the centre of $L$ belongs to class $\Omega$; 3)
$\frak M$ belongs to class $\Omega$; 4) the centre of $\frak M$
belongs to class $\Omega$.

There  exists a commutative Moufang loop (CML) with trivial centre
[5]. Then for the described  CML with finiteness conditions it is
reasonable to use the notion of centralizer, more general that the
notion of centre. This paper generalizes the aforementioned result
for $ZA$-loops. It is proved that for a CML $L$ with
multiplication group $\frak M$ the following statements are
equivalent: 1) the CML $L$ belongs to class $\Omega$; 2) the
centralizer of some finitely generated subloop of $L$ belongs to
class $\Omega$; 3) the group $\frak M$ belongs to class $\Omega$;
4) the centralizer of some finitely generated subgroup of group
$\frak M$ belongs to class $\Omega$.

We denote that the equity of statements 3) and 4) is false for
arbitrary group . In [6] it is constructed an example of locally
nilpotent group of infinite rank in which the centralizer is an
infinite cyclic subgroup and has a finite rank. The centre of this
group is different from unity.

Let $G$ be a group. The commutator $[a,b]$ of the elements $a, b
\in L$ is defined by equality $[a,b] = a^{-1}b^{-1}ab =
a^{-1}a^b$. The identity
$$[xy,zt] = [x,t]^y[y,t][x,z]^{yt}[y,z]^t \eqno{(1)}$$
holds in $G$.

Let $M$ be a subset and $H$ be a subgroup of group $G$. The
subgroup $C_H(M) = \{x \in H \vert [x,y] = 1 \quad \forall y \in
M\}$ is called \textit{centralizer} of subset $M$ into subgroup
$H$. If $M = H = G$, then the normal subgroup $C_G(G)$ is called
\textit{centre} of group $G$.  We will denote it by $C(G)$.

Let us bring some notions and results on the theory of the
\textit{commutative Moufang loops} (abbreviated \textit{CML}s)
from [5], which are characterized by the identity $x^2\cdot yz =
xy\cdot xz.$

The \textit{multiplicative group} $\frak M(L)$ of a CML $L$ is the
group generated by all the \textit{translations} $L(x)$, where
$L(x)y = xy$. The \textit{associator} $(a,b,c)$ of the elements
$a, b, c$ in CML $L$ is defined by the equality $ab\cdot c = (a
\cdot bc)(a,b,c)$. The identities:
$$(xy,u,v) = (x,u,v)((x,u,v),x,y)(y,u,v)((y,u,v),y,x),
\eqno{(2)}$$
$$(x,y,z) = (y^{-1},x,z) = (y,x,z)^{-1} = (y,z,x)\eqno{(3)}$$
hold  in  CML $L$.

The \textit{centre} $Z(L)$ of a CML $L$ is the normal subloop
$Z(L) = \{x \in L \vert (x,y,z) = 1 \quad \forall y,z \in L\}$.
The \textit{upper central series} of the CML $L$ is the series

$$1 = Z_0 \subseteq Z_1 \subseteq Z_2 \subseteq \ldots \subseteq
Z_{\alpha} \subseteq \ldots$$ of the normal subloops of the CML
$L$, satisfying the conditions: 1) $Z_{\alpha} = \sum_{\beta <
\alpha} Z_{\beta}$ for the limit ordinal and 2) $Z_{\alpha +
1}/Z{\alpha} = Z(L/Z_{\alpha})$ for any $\alpha$. If the CML
possesses a central series, then this loop is called
\textit{$ZA$-loop}. If the upper central series of the $ZA$-loop
has a finite length, then the loop is called \textit{centrally
nilpotent}. The least of such a length is called the
\textit{class} of the central nilpotentce.

By analogy it is defined the notion of \textit{$ZA$-group} with
the help of the notion of centre of group.
\smallskip\\
\textbf{Lemma 1 (Bruck-Slaby Theorem)}. \textit{Let $n$ be a
positive integer, $n \geq 3$. Then every commutative Moufang loop
$L$ which can be generated by $n$ elements is centrally nilpotent
of class at most $n - 1$.}
\smallskip\\
\textbf{Lemma $1^0$.} \textit{The multiplication group of any CML
is locally nilpotent.}
\smallskip\\
\textbf{Lemma 2}. \textit{For any CML $L$ with centre $Z(L)$ the
quotient loop $L/Z(L)$ is locally finite $3$-loop of exponent $3$
and is finite if $L$ is finitely generated}.
\smallskip\\
\textbf{Lemma $2^0$}. \textit{Let $L$ be a CML with multiplication
group $\frak M$ and let $C(\frak M)$ be the centre of $\frak M$.
Then the quotient group $\frak M/C(\frak M)$ is locally finite
$3$-group and is finite if $L$ is finitely generated.}

The following concept is the natural generalization of the concept
of centre. Let $M$ be a subset and $H$ be a subloop of CML $L$.
The set $Z_H(M) = \{x \in H \vert x\cdot yz = xy\cdot z \quad
\forall y, z \in M\}$ is called \textit{centralizer} of subset $M$
into subloop $H$. $Z_H(M)$ is a subloop of $L$ [7].
\smallskip\\
\textbf{Lemma 3.} \textit{A centrally nilpotent CML $L$ belongs to
class $\Omega$ if and only if the centralizer of some of its
finitely generated subloop $H$ also belongs to class $\Omega$.}
\smallskip\\
\textbf{Proof.} The necessity of lemma is obvious. To prove the
sufficiently, it is enough  to achieve  by induction by class of
central nilpotence of CML $L$. We will prove this only for the
case of rank finiteness, as the proof of other cases of class
$\Omega$ follows the same pattern.

Let the subloop $H$ be generated by elements $a_1, a_2, \ldots,
a_n$ and let the centralizer $Z_L(H)$ have a finite rank. We will
suppose that the CML $L$ is non-associative, as for abelian groups
(centrally nilpotent CML of class $k = 1$) the statement holds.

Let $Z$ be the centre of CML $L$ and let $k$ be the class of
central nilpotence. Obviously, $Z \subseteq Z_L(H)$, hence the
rank of $Z$ is finite. As the quotient loop $L/Z$ is centrally
nilpotent of class $k - 1$, then by inductive assertion the rank
of $L/Z$ will be finite if the centralizer $D/Z$ of image $HZ/Z$
of subloop $H$ into quotient loop $L/Z$ has a finite rank. We will
prove this.

Let $x, y \in D$, and let $A_i = \{a_{i_1}, a_{i_2}\}, 1 = 1, 2,
\ldots, t$, be an arbitrary fixed pair of elements $a_{i_1},
a_{i_2} \in A$. We have $(x,a_{i_1},a_{i_2}) \in Z$, then by (1)
$(xy,a_{i_1},a_{i_2}) = (x,a_{i_1},a_{i_2})(y,a_{i_1},a_{i_2})$.
This equality shows that the mapping $x \rightarrow
(x,a_{i_1},a_{i_2})$ is a homomorphism of $D$ into $Z$. For each
$A_i$ we consider the homomorphisms $\varphi_i(x) =
(x,a_{i_1},a_{i_2})$, $x \in D$. Obviously, $\text{ker}\varphi_i =
Z_D(A_i)$. Then the quotient loop is an abelian group of finite
rank. The direct product $\prod^t_{i=1}D/Z_D(A_i)$ also has a
finite rank. By (1), (2) it is easy to see that
$\bigcap_{i=1}^tZ_D(A_i) = Z_D(A)$.

Analogously to Remak Theorem for groups [8] it may be proved that
the quoti\-ent loop $D/Z_D(H)$ is isomorphic to subloop of direct
product \break $\prod^t_{i=1}D/Z_D(A_i)$. Then $D/Z_D(H)$ has a
finite rank. As $Z_D(H) \subseteq Z_L(H)$ then $Z_D(H)$ also has a
finite rank. Hence  the CML $D$ also has a finite rank.

Consequently, $L/Z$ is a CML of finite rank. We mentioned above
that the centre $Z$ has a finite rank. Then  CML $L$ also  has a
finite rank. This completes the proof of Lemma 3.
\smallskip\\
\textbf{Lemma $3^0$} [6]. \textit{A  nilpotent group  belongs to
class $\Omega$ if and only if the centralizer of some its finitely
generated subloop belongs to class $\Omega$.}
\smallskip\\
 \textbf{Lemma 4.} \textit{Let $H$ be a finitely generated subloop
 of CML $L$. If the centralizer $Z_L(H)$  belongs to class $\Omega$ then
 the centralizer \break $Z_{L/Z(L)}(Z(L)H/Z(L))$ belongs to class $\Omega$.}
\smallskip\\
\textbf{Proof.} In [9, 10] it is proved that for CML the condition
of finitely generation  and maximum condition for subloops are
equivalent, and in [3] is proved that these conditions are
equivalent with the maximum condition for associative subloops.
Further, in [2] it is proved that for a CML the minimum condition
for subloops and the minimum condition for associative subloops
are equivalents, and for $p$-loops this conditions are equivalents
with condition of finiteness rank [4].

Now, let $Z(L)$ be the centre of CML $L$, let $\overline L =
L/Z(L)$, $\overline H = \break Z(L)H/Z(L)$, and we suppose that
the centralizer $Z_{\overline L}(\overline H)$ does  not belong to
class $S$. By Lemma 2 the quotient loop $L/Z(L)$ satisfies the
identity $x^3 = 1$. Then by the aforementioned,  $Z_{\overline
L}(\overline H)$ contains an infinite elementary abelian $3$-group
$\overline B$, which decomposes into a direct product of cyclic
groups of order $3$. Let $A/Z(L) = \overline A = \overline A_1
\times \overline A_2 \times \ldots \times \overline A_i \times
\ldots $ be the maximal subgroup of $\overline B$ regarding to
property $\overline A_i \nsubseteq \overline H$, and let
$\overline A_i = <\overline a_i>$. We denote by $M(\overline a_i)$
a maximal subloop of CML $\overline R = <\overline A, \overline
H>$ such that $\overline a_i \notin M(\overline a_i)$. As the
element $\overline a_i$ has order $3$ then $\overline A_i \cap
M(\overline a_i) = \overline 1$. Every maximal subloop of CML is
normal in this CML [10]. Let $I(\overline R)$ be the inner mapping
group of CML $\overline R$. Every inner mapping of CML is an
automorphism of this CML [5]. Then $I(\overline R)\overline A_i =
\overline A_i$. Hence the subloop $\overline A_i$ is normal in
$\overline R$. In [12] it is proved that if an element of order
$3$ of CML generates a normal subloop then this element belongs to
the centre of this CML. Hence $\overline A_i \subseteq Z(\overline
R)$. $\overline A \subseteq Z(\overline R)$. From here it follows
that $\overline R = \overline A\overline H$. But $\overline A \cap
\overline H = \overline 1$. Then from $I(\overline R)\overline A =
\overline A$ it follows that $I(\overline R)\overline H =
\overline H$, i.e. the subloop $\overline H$ is normal in
$\overline R$. Consequently, $\overline R = \overline A \times
\overline H$.

The subloop $H$ is finitely generated. Then by Lemma 1 its is
centrally nilpotent. The subloop $\overline H$ is also centrally
nilpotent. Then and subloop $\overline R$ is centrally nilpotent.
As $\overline B \subseteq \overline R$ and $\overline B$ do not
belong to class $S$ then  $\overline R$ also does not belong to
class $S$. The inverse image of $\overline R$ under the
homomorphism $L \rightarrow L/Z(L)$ is $AH$. This CML is centrally
nilpotent and does not belong to class $S$. Then by Lemma 3  the
centralizer $Z_{AH}(H)$ does not belong to class $S$. We get a
contradiction, as $Z_{AH}(H) \subseteq Z_L(H)$ and $Z_L(H)$ belong
to class $S$. Consequently, the centralizer $Z_{\overline
L}(\overline H)$ belongs to class $S$.  This completes the proof
of Lemma 4.
\smallskip\\
\textbf{Lemma $4^0$.} \textit{Let $\frak N$ be a finitely
 generated subgroup of multiplication group $\frak M$  of CML $L$.
 If the centralizer $C_{\frak M}(\frak N)$  belongs to class $\Omega$
 then   the centralizer $C_{\frak M/C(\frak M)}(C(\frak M)\frak N/C(\frak
 M)$ also  belongs  to class $\Omega$.}
\smallskip\\
\textbf{Proof.} In [3] it is proved that the finiteness of
generators, the maximum condition for subgroups and maximum
condition for abelian subgroups are equivalent for group $\frak
M$. Further, in [2] it is proved that for $\frak M$ the minimum
condition for subgroups and the minimum condition for abelian
subgroups are equivalent, and for $p$-groups these conditions are
equivalent with the condi\-tion of finiteness rank [4, 11].

We suppose that the centralizer $\frak B/C(\frak M) = C_{\frak
M/C(\frak M)}(C(\frak M)\frak N/C(\frak M))$ does not belong to
class $\Omega$. Then by the first paragraph $\frak B/C(\frak M)$
contains an abelian subgroup $\frak A/C(\frak M)$ which does not
belongs to class $\Omega$. From definition of centrali\-zer it is
follows that $[\frak a,\frak n] = 1$ for $\frak a \in \frak A,
\frak n \in \frak N$. Using this, it is easy to see that the
product $\frak A\frak N$ is a group. Moreover, $\frak A$ is
abelian, $\frak N$ is nilpotent, then also  using  (1) we prove
that $\frak A\frak N$ is a nilpotent group. $\frak A$ does not
belongs to class $\Omega$ then and $\frak A\frak N$ does not
belongs to class $\Omega$. Then by Lemma $3^0$ the centralizer
$C_{\frak A\frak N}(\frak N)$ does not belongs to class $\Omega$.
$C_{\frak M}(\frak N)$ also does not belongs to class $\Omega$ as
$C_{\frak A\frak N}(\frak N) \subseteq C_{\frak M}(\frak N)$. We
get a contradiction. Consequently, $C_{\frak M/C(\frak M)}(C(\frak
M\frak N/C(\frak M)$ belongs
 to class $\Omega$, as required.
\smallskip\\
\textbf{Lemma 5.} \textit{Let $H$ be a finitely generated subloop
 of CML $L$. If the centrali\-zer $Z_L(H)$ belongs to class $\Omega$ then
 the centre $Z(L)$ if CML $L$ is different from unity.}
 \smallskip\\
\textbf{Proof.} If $a \in L$ is an element of infinite order then
 by Lemma 2 $1 \neq a^3 \in Z(L)$. Then we will suppose that $L$
 is a  periodic CML. In this cases $L$ decomposes into a direct product
 of its maximal
$p$-subloops $L_p$, in addition  $L_p$ belongs to the centre
$Z(L)$ under $p \neq 3$. Hence to prove  Lemma 3 it is sufficient
to suppose that $L$ is a $3$-loop. By Lemma 1  every finitely
generated CML is  centrally nilpotent, then we will suppose that
CML $L$ is infinite.

We remind that a system $\{G_{\alpha}\}$ ($\alpha \in I$) of
subloops of loop $G$ is a \textit{local system} if the union
$\bigcup_{\alpha \in I}G_{\alpha}$ coincides with $G$ and every
two members of this system are contained in a certain third member
of this system. Using the definition of local system it is easy to
prove the statement: if $\{G_{\alpha}\}$ ($\alpha \in I$) is some
local system of loop $G$ and $I = I_1 \bigcup I_2 \bigcup \ldots
\bigcup I_k$ is a certain partition of the set of indices $I$ into
a finite number of subsets, $I_j$, $j = 1, 2, \ldots, k$, then at
least on subset $I_j$ corresponds to the set of subloops
$\{G_{\beta}\}$, $\beta \in I_j$, which will also be a local
system for  loop $G$.

Let now $\{L_{\alpha}\}$, $\alpha \in I$, be the local system of
all finitely generated subloops of CML $L$, which contain the
subloop $H$. By Lemma 1 $Z(L_{\alpha}) \neq \{1\}$. For each
$\alpha \in I$ we fixed an arbitrary non-unitary element
$a_{\alpha}$ from centre $Z(L_{\alpha})$ and let $K$ be the
subloop of $L$ generated by all $a_{\alpha}$, $\alpha \in I$. From
$K \subseteq Z_L(H)$ it follows that  $K$ belongs to class
$\Omega$. We suppose that $3$-subloop $K$  has a finite rank. Then
 by [4], $K$ satisfies the minimum condition for its subloops and
by [2]  $K = R \times T$, where $R \subseteq Z(L)$ and $T$ is a
finite subloop. Further, any periodic CML is locally finite [5].
Hence to prove  Lemma 3 it is sufficient to consider that $K$ is a
finite subloop.

Further, let us decompose the set of indices $I$ into a finite
number of subsets $I = I_1 \bigcup I_2 \bigcup \ldots \bigcup I_k$
by rule: $\beta, \gamma \in I_j$ if and only if $a_{\beta} =
a_{\gamma}$. According to the aforementioned statement we have
received that at least one of the subsets $I_j$ (e.g. $ I_1$)
corresponds to the subset of subloops $L_{\alpha}$, $\alpha \in
I_1$, which will be a local system for  CML $L$. Next, let us fix
index  $\alpha \in I_1$, and consider the set of indices  $S
\subseteq I_1$, such that $L_{\alpha} \subseteq L_{\beta}$, $\beta
\in S$. We notice that the set  $S$ corresponds to the set of
subloops $\{L_{\beta}\}$, $\beta \in S$, which will be also a
local system for  CML $L$. Let us denote the value of
corresponding members by  $b$, $b = a_{\alpha} = a_{\beta} =
\ldots$. Then $b \in Z(L_{\beta})$ for all $\beta \in S$ and,
consequently, $b \in Z(L)$. This completes the proof of Lemma 5.
\smallskip\\
\textbf{Lemma $5^0$.} \textit{Let  $\frak N$ be a finitely
generated subgroup of multiplication group $\frak M$
 of CML $L$. If the centrali\-zer
 $C_{\frak M}(\frak N$) belongs to class $\Omega$ then
 the centre  $C(\frak M)$
 of group $\frak M$ is different from unity.}
 \smallskip\\
 \textbf{Proof.} We will prove Lemma $5^0$, using the same pattern
 as for Lemma 5.  If  $a \in \frak M$ is
 an element of infinite order then
 by Lemma  $2^0$  $a^{3^k} \in C(\frak M)$ for some integer $k$.
 Then we will suppose that $\frak M$
 is a  periodic group. In this case $\frak M$
  decomposes into a direct product
 of its maximal
$p$-subgroups  $\frak M_p$, in addition  $\frak M_p$ belongs to
the centre  $C(\frak M)$ under $p \neq 3$. Hence to prove Lemma
$5^0$ it is sufficient to suppose that  $\frak M$ is a $3$-group.
By Lemma $1^0$ a finitely generated subgroup of multiplication
group of CML is nilpotent, then we will suppose that $\frak M$ is
infinite.

Let  $\{\frak M_{\alpha}\}$, $\alpha \in I$, be a local system of
all finitely generated subgroups of group $\frak M$, which contain
the subgroup $\frak N$. By Lemma $1^0$ $C(\frak M_{\alpha}) \neq
\{1\}$. For each $\alpha \in I$ we fixed an arbitrary non-unitary
element $a_{\alpha}$ from centre $C(\frak M_{\alpha})$ and let
$\frak K$ be the sugroup of $\frak M$ generated by all
$a_{\alpha}$, $\alpha \in I$. From $\frak K \subseteq C_{\frak
M}(\frak N)$ it follows that $\frak K$ belong to class $\Omega$.
We suppose that $3$-subgroup $\frak K$ has a finite rank. Then
 by [4, 11] $\frak K$ satisfies the minimum condition for its subgroups and
by [2]  $\frak K = \frak R \times \frak T$, where $R \subseteq
C(\frak M)$ and $\frak T$ is a finite subloop. Any periodic
multiplication group of CML is locally finite [5]. Then to prove
Lemma $5^0$ it is sufficient to consider that $\frak K$ is a
finite subgroup.

Further, let us decompose the set of indices $I$ into a finite
number of subsets $I = I_1 \bigcup I_2 \bigcup \ldots \bigcup I_k$
by the rule: $\beta, \gamma \in I_j$ if and only if $a_{\beta} =
a_{\gamma}$. According to the aforementioned statement we have
received that at least one of the subsets $I_j$ (e.g. $ I_1$)
corresponds to the subset of subgroups $\frak M_{\alpha}$, $\alpha
\in I_1$, which will be a local system for group $\frak M$. Next,
let us fix index  $\alpha \in I_1$, and consider the set of
indices $S \subseteq I_1$, such that $\frak M_{\alpha} \subseteq
\frak M_{\beta}$, $\beta \in S$. We notice that the set  $S$
corresponds to the set of subgroups $\{\frak M_{\beta}\}$, $\beta
\in S$, which will be also a local system for group $\frak M$. Let
us denote the value of corresponding members by  $b$, $b =
a_{\alpha} = a_{\beta} = \ldots$. Then $b \in C(\frak M_{\beta})$
for all $\beta \in S$ and, consequently, $b \in C(\frak M)$. This
completes the proof of Lemma $5^0$.
\smallskip\\
\textbf{Theorem.} \textit{For a CML $L$ with multiplication group
$\frak M$ the following statements are equivalent:}

\textit{1) the CML $L$ belongs to class $\Omega$;}

{2) the centralizer of some finitely generated subloop $H$ of $L$
belongs to class $\Omega$;}

\textit{3) the group $\frak M$ belongs to class $\Omega$;}

\textit{4) the centralizer of some finitely generated subgroup
$\frak N$ of group $\frak M$ belongs to class $\Omega$.}
\smallskip\\
\textbf{Proof.} $2) \Rightarrow 1)$. Let the centralizer $Z_L(H)$
belong to class $\Omega$. We denote $\overline L = L/Z(L)$,
$\overline H = HZ(L)/Z(L)$. Any periodic CML is locally finite
[5]. Then from Lemma 2 it follows that the subloop $\overline H$
is finite. From Lemmas 4, 5 it follows that the upper central
series of CML $\overline L$ has a form $\overline 1 \subset
Z_1(\overline L) \subset \ldots \subset Z_k(\overline L) \subset
\ldots $, and for a natural number $n$ $Z_n(\overline L) \neq
Z_{n+1}(\overline L)$ if $Z_n(\overline L) \neq \overline L$. As
$\overline H$ is finite, then for some $k$ $\overline H \nsubseteq
Z_{k-1}(\overline L)$, but $\overline H \subseteq Z_k(\overline
L)$. Then $Z_{\tilde{L}}(\tilde{H}) = \tilde{L}$, where $\tilde{L}
= \overline{L}/\overline{Z_k}$, $\tilde{H} =
\overline{H}\overline{Z_k}/\overline{Z_k}$. By Lemma 2 $\tilde{L}$
is a $3$-loop and by Lemma 4 $\tilde{L}$ belongs to class
$\Omega$. In [4] it is proved that the minimum condition for
subloops and the condition of finiteness rank are equivalent for
CML $\tilde{L}$. In this casse $\tilde{L} = \tilde{R} \times
\tilde{T}$, where $\tilde{R} \subseteq Z(\tilde{L})$, and
$\tilde{T}$ is a finite CML which by Lemma 1 is  centrally
nilpotent. Then  CML $\tilde{L} = \overline{L}/\overline{Z_k} =
(L/Z)/(Z_k/Z) \cong L/Z_k$ is also centrally nilpotent. From
central nilpotence of $L/Z_k$ it follows the central nilpotence of
$L$. By Lemma 3 $L$ belongs to class $\Omega$. Consequently, the
implication $2) \Rightarrow 1)$ hold.

4) $\Rightarrow 3)$. Let the centralizer $C_\frak M(\frak N)$
belong to class $\Omega$. We denote $\overline{\frak M} = \frak
M/C(\frak M)$, $\overline{\frak N} = \frak NC(\frak M)/C(\frak
M)$. Any periodic multiplication group of CML is locally finite
[5]. Then from Lemma $2^0$ it follows that the subgroup
$\overline{\frak N}$ is finite. From Lemmas $4^0$, $5^0$ it
follows that the upper central series of group $\overline{\frak
M}$ has a form $\overline{1} \subset C_1(\overline{\frak M})
\subset \ldots \subset C_k(\overline{\frak M}) \subset \ldots $,
and for a natural number $n$ $C_n(\overline{\frak M}) \neq
C_{n+1}(\overline{\frak M})$ if $C_n(\overline{\frak M}) \neq
\overline{\frak M}$. As $\overline{\frak N}$ is finite, then for
some $k$ $\overline{\frak N} \nsubseteq C_{k-1}(\overline{\frak
M})$, but $\overline{\frak N} \subseteq C_k(\overline{\frak M})$.

Then $C_{\tilde{\frak M}}(\tilde{\frak N}) = \tilde{\frak M}$,
where $\tilde{\frak M} = \overline{\frak M}/\overline{C_k}$,
$\tilde{\frak N} = \overline{\frak
N}\overline{C_k}/\overline{C_k}$. By Lemma $2^0$ $\tilde{\frak M}$
is a $3$-loop and by Lemma $4^0$ $\tilde{\frak M}$ belongs to
class $\Omega$. In [4] it is proved that the minimum condition for
subgroups and the condition of finiteness rank are equivalent for
group $\tilde{\frak M}$. In this case $\tilde{\frak M} =
\tilde{\frak R} \times \tilde{\frak T}$, where $\tilde{\frak R}
\subseteq C(\tilde{\frak M})$, and $\tilde{\frak T}$ is a finite
group which by Lemma $1^0$ is  nilpotent. Then  group
$\tilde{\frak M} = \overline{\frak M}/\overline{C_k} = (\frak
M/C)/(C_k/C) \cong \frak M/C_k$ is also  nilpotent. From
nilpotence of $\frak M/C_k$ it follows the  nilpotence of $\frak
M$. By Lemma $3^0$ $\frak M$ belongs to class $\Omega$.
Consequently, the implication $4) \Rightarrow 3)$ hold.

Hence $C_{k+1}(\overline{\frak M}) = \overline{\frak M}$. From
here it follows that the group $\overline{\frak M}$ is  nilpotent.
Then and the group $\frak M$ is nilpotent and by Lemma $3^0$ $L$
belongs to class $\Omega$.

The implications $1) \Rightarrow 2)$, $3) \Rightarrow 4)$ are
obvious.

The statements: $L$ is finitely generated and the maximum
conditions holds in $L$ are equivalent for any CML $L$ [9, 10].
Then the implications $1) \Leftrightarrow 3)$ are proved in [2 --
4]. This completes the proof of the Theorem.

Tiraspol State University (Moldova)

E-mail: aliona2010@yahoo.md; sandumn@yahoo.com

\end{document}